\documentclass{amsart}
\usepackage{amssymb}
\newtheorem{theorem}{Theorem}[section]
\newtheorem{lemma}[theorem]{Lemma}
\newtheorem{cor}[theorem]{Corollary}
\newtheorem{proposition}[theorem]{Proposition}
\newtheorem*{thmA}{Theorem A}
\newtheorem*{thmB}{Theorem B}
\newtheorem*{3sbgrp}{Three Subgroup Lemma}

\theoremstyle{definition}

\theoremstyle{remark}

\begin{document}

\title{Generalizing a theorem of P. Hall on finite-by-nilpotent groups}
\author{Gustavo A. Fern\'andez-Alcober}
\address{Matematika Saila, Euskal Herriko Unibertsitatea,
48080 Bilbao, Spain }
\email{gustavo.fernandez@ehu.es}
\author{Marta Morigi}
\address{Dipartimento di Matematica, Universit\`a di Bologna, Piazza
di Porta San Donato 5, 40127 Bologna, Italy}
\email{mmorigi@dm.unibo.it}

\thanks{The first author is supported by the Spanish Ministry of Science and Education,
grant MTM2004-04665, partly with FEDER funds, and by the University of the Basque
Country, grant UPV05/99. The second author is partially supported by
MIUR (Project ``Teoria dei Gruppi e applicazioni'') and thanks the
University of the Basque Country for the hospitality.}



\begin{abstract} Let $\gamma_i(G)$ and $Z_i(G)$ denote the
  $i$-th terms of  the lower and upper central series of a group $G$,
  respectively.  In
  \cite{H} P. Hall showed
  that if $\gamma_{i+1}(G)$ is finite then the index $|G:Z_{2i}(G)|$
  is finite. We
  prove that the same result holds under the weaker hypothesis that
  $|\gamma_{i+1}(G):\gamma_{i+1}(G)\cap Z_i(G)|$ is finite.
\end{abstract}
\maketitle

\section{Introduction}

If $G$ is an arbitrary group, a classical theorem of Schur asserts that
if the center $Z(G)$  has finite index in $G$ then derived
subgroup $G^\prime$ is finite. This was later generalized by
Baer (see 14.5.1 in \cite{R}) to any term of the lower central series;
namely, if $|G:Z_i(G)|$ is finite then $\gamma_{i+1}(G)$ is finite.
The converse does not hold in general, anyway in \cite{H} P. Hall proved
that if $\gamma_{i+1}(G)$ is finite then  $|G:Z_{2i}(G)|$ is finite.
For the case $i=1$ a stronger result is known to hold: actually, if
$|G^\prime:G^\prime\cap  Z(G)|$ is finite then $|G:Z_{2}(G)|$ is also
finite. This result was obtained independently  by the first author and
Moret\'o  (see Theorem E of \cite{FM}) and by Podoski and Szegedy
in \cite{PS}. In this paper we show
that this last property can be extended to an arbitrary value of
$i$. More precisely, the following is true.

\begin{thmA}
Let $G$ be  a group such that
$|\gamma_{i+1}(G):\gamma_{i+1}(G)\cap Z_i(G)|$ is
finite. Then $|G:Z_{2i}(G)|$ is finite.
\end{thmA}

From the proof of the theorem it can be checked that $|G:Z_{2i}(G)|$ is
bounded in terms of $|\gamma_{i+1}(G):\gamma_{i+1}(G)\cap
Z_i(G)|$, the ultimate reason being that Lemma \ref{fg} and Theorem
\ref{baer} below can be stated in a quantitative version. However, we have made no
attempt at giving a sharp bound. We will mention here that in the case
$i=1$, the existence of such a bound  was proved by Isaacs in \cite{I} when the group $G$
is finite, and then an explicit bound was given in \cite{PS} for an
arbitrary group $G$.

Related to Theorem A, the following two questions arise naturally:
\begin{enumerate}
\item
To what extent is Theorem A best possible?
If the weaker condition that $|\gamma_{i+1}(G):\gamma_{i+1}(G)\cap
Z_{i+1}(G)|$ is finite holds, can we conclude that $|G:Z_{2i}(G)|$ is finite?
\item
In the case that $\gamma_{i+1}(G)$ is finite, if $G$ is also finitely generated,
then a stronger result holds, namely $|G:Z_i(G)|$ is finite.
Is this true also under the hypothesis of Theorem A?
Does it follow at least that $|G:Z_j(G)|$ is finite for some $j$ smaller than $2i$?
\end{enumerate}
The answer to both these questions is negative.
To see this, for arbitrary $c$, consider a finitely generated nilpotent group $G$ of class
$c$ in which the upper and lower central series coincide and such that $|G:Z_{c-1}(G)|$ is infinite.
For example, one can take the semidirect product $G=B\ltimes A$, where
$B=\langle b\rangle$ is an infinite cyclic group,
$A$ is the free abelian group on free generators $a_1,\ldots,a_c$, and $b$ acts on $A$ by $a_i^b=a_ia_{i+1}$
for $1\le i\le c-1$ and $a_c^b=a_c$.
Now if $i\ge 1$ is any fixed integer, we get counterexamples to the first and the second questions
by choosing $c=2i+1$ and $c=2i$, respectively.

Finally, we observe that combining our result with Baer's theorem it follows
that if $|\gamma_{i+1}(G):\gamma_{i+1}(G)\cap Z_i(G)|$ is
finite then $\gamma_{2i+1}(G)$ is also finite.
Actually, one of the key arguments in our proof of Theorem A is the following
generalization of this fact, which might be interesting in its own right.

\begin{thmB}
Let $G$ be a group such that
$|\gamma_{s}(G):\gamma_{s}(G)\cap Z_t(G)|$ is
finite for some $s,t$. Then $|\gamma_{s+j}(G):\gamma_{s+j}(G)\cap
Z_{t-j}(G)|$ is finite
for every $j$ such that $0\le j\le t$. In particular,
$\gamma_{s+t}(G)$ is finite.
\end{thmB}

\section{The results}

The notation we use is standard. Moreover, following the book
\cite{R}, if $A$ and $B$ are subgroups of a group
$G$ and $n$ is a natural number, we define recursively:
$$[A,\,_1B]=[A,B],\qquad\qquad [A,\,_nB]=[A,\,_{n-1}B,B].$$

Throughout the paper, we will repeatedly use the following well-known
result (see for instance 5.1.10 in \cite{R}).

\begin{3sbgrp}
Let $H,K,L$ be
  subgroups of a
  group $G$. If two of the commutator subgroups
  $[H,K,L],[K,L,H],[L,H,K]$ are contained in a normal subgroup of $G$,
  then so is the third.
\end{3sbgrp}

 Another result which will be often used in our proofs is stated for
 convenience in the following lemma, whose proof is elementary. Most
 of the times, we will apply it modulo a normal subgroup.

 \begin{lemma}\label{fg}
Let $H,K$ be subgroups of a group $G$. If $[H,K]$ is finite and $H$ is
finitely generated, then the centralizer $C_K(H)$  has
finite index in $K$.
\end{lemma}

We will also need the following result of Baer (see for instance 14.5.2
in \cite{R}).

\begin{theorem}\label{baer} Let
$M\le H$ and $N\le K$ be normal subgroups of a group $G$ such that
$|H:M|$  and $|K:N|$ are finite, and
$[H,N]=1=[K,M]$. Then $[H,K]$ is finite.
\end{theorem}

The key step in the proof of our main theorem is in the following
proposition.

\begin{proposition}\label{gf}
Let $G$ be a group and let $s\ge1$ be an integer  such that
$|\gamma_{s}(G):\gamma_{s}(G)\cap Z(G)|$ is finite. Then
$C_G\left(\gamma_{s}(G)\right)$ has finite index in $G$ and
$\gamma_{s+1}(G)$ is finite.
\end{proposition}

\begin{proof} Let $Z= \gamma_{s}(G)\cap Z(G)$. As
$\gamma_{s}(G)/Z$ is
finite, there exists a finitely generated subgroup $U$ of $G$ such
that $\gamma_{s}(G)=\gamma_{s}(U)Z$.
By applying P. Hall's theorem to the quotient group
$G/Z(G)$, we obtain that
$|G:Z_{2s-1}(G)|$ is finite. By the theorem of Baer mentioned in the introduction, it follows
that $\gamma_{2s}(G)$ is finite. Since
$\gamma_k(U)/\gamma_{k+1}(U)$ is finitely generated for every
$k=1,\ldots,2s-1$, we conclude that all terms of the lower central series
of $U$ are finitely generated.

We are going to prove that, for every $j=1,\ldots,s$, there exists a
subgroup $H_j$ of finite index in $\gamma_j(G)$ such that
$[H_j,\gamma_{s-j+1}(U)]=1$. Then $[H_1,\gamma_{s}(G)]=[H_1,\gamma_s(U)Z]=
[H_1,\gamma_s(U)]=1$, which proves that  $|G:C_G(\gamma_{s}(G))|$ is finite.

We prove the existence of $H_j$ by reverse induction on $j$.
For $j=s$, we take $H_s=Z$.
Suppose now that we already have $H_{j+1}$ of finite index in $\gamma_{j+1}(G)$
such that $[H_{j+1},\gamma_{s-j}(U)]=1$, and let us see how to construct the
subgroup $H_j$. Let
$K_j=C_{\gamma_j(G)} (\gamma_{s-j}(U)Z/Z)$.
Since $\gamma_{s-j}(U)Z/Z$ is
finitely generated and $[\gamma_j(G),\gamma_{s-j}(U)]Z/ Z  \le
\gamma_{s}(G)/Z$ is finite, it follows from Lemma \ref{fg} that $K_{j}$ has
finite index in $\gamma_j(G)$. Also
\begin{equation}
\label{Kj}
[K_j,\gamma_{s-j}(U),U]\le [Z,U]=1.
\end{equation}

Let now $D_{j+1}=C_{\gamma_{j+1}(G)}(\gamma_{s-j}(U))$.
Since $[H_{j+1},\gamma_{s-j}(U)]=1$, we have $H_{j+1}\le D_{j+1}$
and consequently $|\gamma_{j+1}(G):D_{j+1}|$ is finite.
Consider $T_{j}=\gamma_{j}(G)U$. We claim that $D_{j+1}$ is normal in $T_j$.
On the one hand, for every $u$ in $U$ we have
$[D_{j+1}^u,\gamma_{s-j}(U)]=[D_{j+1},\gamma_{s-j}(U)]^u=1$, so
that $U$ normalizes $D_{j+1}$.
On the other hand, we have
$[\gamma_{j}(G),\gamma_{s-j}(U),D_{j+1}]\le[\gamma_{s}(G),D_{j+1}]
=[\gamma_{s}(U)Z,D_{j+1}]=1$ and
$[\gamma_{s-j}(U),D_{j+1},\gamma_{j}(G)]=1$ by the definition of
$D_{j+1}$. By the Three Subgroup Lemma, $[D_{j+1},\gamma_{j}(G),
\gamma_{s-j}(U)]=1$ and consequently also $\gamma_{j}(G)$ normalizes
$D_{j+1}$.

Now we work in the quotient group
$T_j/D_{j+1}$. Since
$[\gamma_{j}(G),U]D_{j+1}/D_{j+1}\le\gamma_{j+1}(G)/D_{j+1}$ is finite and
$U$ is finitely generated, it follows that the centralizer
$L_j$ of $UD_{j+1}/D_{j+1}$ in $\gamma_{j}(G)$ has finite index in
$\gamma_{j}(G)$. Observe that
\begin{equation}\label{Lj}[L_j,U,\gamma_{s-j}(U)]\le[D_{j+1},\gamma_{s-j}(U)]=1.
\end{equation}

Finally, let $H_j=K_j\cap L_j$. Then $|\gamma_j(G):H_j|$ is
finite. Moreover, using (\ref{Kj}) and (\ref{Lj}) and the Three
Subgroup Lemma we obtain that $[\gamma_{s-j}(U),U,H_j]=1$, that is
$[\gamma_{s-j+1}(U),H_j]=1$, as desired.

Now in order to prove that $\gamma_{s+1}(G)$ is finite we apply Theorem
\ref{baer} with $M=C_G(\gamma_s(G))$, $H=G$, $N=\gamma_s(G)\cap Z(G)$ and
$K=\gamma_s(G)$. It follows that $[G,\gamma_s(G)]=\gamma_{s+1}(G)$ is finite.
\end{proof}

Let us remark that if $N$ is a normal
subgroup of a group $G$ and $|N:N\cap Z(G)|$ is finite, it does not follow
that $|G:C_G(N)|$ or $[N,G]$ are finite.
For example, let $p$ be a prime and let $H$ and $N$ be two elementary abelian $p$-groups
with countable bases $\{x_i\}_{i\ge 1}$ and $\{y_j\}_{j\ge 0}$, respectively.
We define an action of $H$ on $N$ so that $x_i$ centralizes all $y_j$ with $j\ge 1$ and $y_0^{x_i}=y_0y_i$.
Then in the semidirect product $G=H\ltimes N$ we have $|N:N\cap Z(G)|=p$ but both $|G:C_G(N)|$ and $[N,G]$ are infinite.

\begin{cor}\label{fin} Let $G$ be a group such that
$|\gamma_{s}(G):\gamma_{s}(G)\cap Z_t(G)|$ is
finite for some $s,t$. Then $|\gamma_{s+j}(G):\gamma_{s+j}(G)\cap
Z_{t-j}(G)|$ is finite
for every $j$ such that $0\le j\le t$. In particular,
$\gamma_{s+t}(G)$ is finite. \end{cor}

\begin{proof}  By induction on $t$ it suffices to prove that
$|\gamma_{s+1}(G):\gamma_{s+1}(G)\cap Z_{t-1}(G)|$ is
finite. This follows immediately by applying Proposition \ref{gf} to
the quotient
group $G/Z_{t-1}(G)$. \end{proof}

The last part of the proof of Theorem A is inspired from Hall's
ideas.
The main role will be played by a subgroup $C$ with
the two properties
that $C$ has finite index in $G$ and $[C,\,_{s-1}G,C]\le Z_{2i-s}(G)$ for
every $s\ge 1$, with the convention that $Z_j(G)=1$ for $j\le 0$. The following technical lemma will ensure
that $C$ has the second property.

\begin{lemma}\label{comm}
 Let $G$ be a group and let
 $C_j$ be the centralizer in $G$ of $\gamma_{i+j}(G)/\gamma_{i+j}(G)
\cap Z_{i-j}(G)$ for $j=1,\ldots, i$. If
$C=\bigcap_{j=1}^{i}\,C_j$, then $[C,\,_{s-1}G,C]\le Z_{2i-s}(G)$ for
every $s\ge 1$.
\end{lemma}

\begin{proof}  Observe that $C$ is normal in $G$ and so is
$[C,\,_{k}G]$ for every $k$. We first prove by induction on $k$ that
\begin{equation}
\label{comm1}
[[C,\,_{k}G],\gamma_{\ell}(G)]\le
  Z_{2i-k-\ell}(G) \quad \mathrm{for\,\, all}\,\,k\ge
  0\,\,\mathrm{and\,\,for\,\,all}\,\,\ell\ge i+1.
\end{equation}
 By definition of
$C$, we have
$[C,\gamma_{\ell}(G)]\le Z_{2i-\ell}(G)$ for all $\ell\ge i+1$ and this
settles the case $k=0$.
Assume now that the statement is true for  $k$.
We have $[\gamma_{\ell}(G),[C,\,_{k}G],G]
\le [Z_{2i-k-\ell}(G),G]\le Z_{2i-k-\ell-1}(G)$.
Also, $[G,\gamma_{\ell}(G),[C,\,_{k}G]]=
[\gamma_{\ell+1}(G),[C,\,_{k}G]]\le Z_{2i-k-\ell-1}(G)$.
So by the Three
Subgroup Lemma, it follows that
\[
[[C,\,_{k+1}G],\gamma_{\ell}(G)]=
[[C,\,_{k}G],G,\gamma_{\ell}(G)]\le Z_{2i-k-\ell-1}(G),
\]
which proves the
statement for $k+1$.

Now, in order to prove the lemma, we need to show that
$[C,\,_{s-1}G,C,\,_{2i-s}G]=1$ for
every $s=1,\ldots, 2i$. We use the formula in the proof of
14.5.4 of \cite{R}, which says that if $M,N$ are normal subgroups of a group
$G$ then
$[[M,N],\,_{n}G]\le\prod_{\ell=0}^n \big[[M,\,_{n-\ell}G],[N,\,_{\ell}G]\big].$
Applying this with $M=[C,\,_{s-1}G],\;N=C$, we have
$$[C,\,_{s-1}G,C,\,_{2i-s}G]\le\prod_{\ell=0}^{2i-s}\big[[C,\,_{2i-\ell-1}G],
[C,\,_{\ell}G]\big].$$
If $0\le \ell \le i-1$ then $2i-\ell\ge i+1$, and since
$[C,\,_{2i-\ell-1}G]\le \gamma_{2i-\ell}(G)$, we have
$[[C,\,_{2i-\ell-1}G], [C,\,_{\ell}G]\big]=1$ by (\ref{comm1}).
If $\ell\ge i$ then we can argue similarly, since $[C,\,_{\ell}G]\le \gamma_{\ell+1}(G)$.
\end{proof}

\begin{theorem}\label{main}
 Let $G$ be  a group such that
$|\gamma_{i+1}(G):\gamma_{i+1}(G)\cap Z_i(G)|$ is
finite. Then  $|G:Z_{2i}(G)|$ is also finite.
\end{theorem}

\begin{proof}
Let $C_j$ be the centralizer in $G$ of $\gamma_{i+j}(G)/(\gamma_{i+j}(G)
\cap Z_{i-j}(G))$ for $j=1,\ldots, i$.
Since $|\gamma_{i+j}(G):\gamma_{i+j}(G)\cap Z_{i-j+1}(G)|$ is finite by Corollary \ref{fin},
we can apply Proposition \ref{gf} to the quotient group $G/Z_{i-j}(G)$,
and it follows that $|G:C_j|$ is finite.
Let $C=\bigcap_{j=1}^{i}\,C_j$, which has also finite index in $G$.

For every $s=1,\ldots,i+1$, put $K_s=[C,\,_{s-1}G]$, which is contained in $\gamma_s(G)$.
We prove by reverse induction on
$s$ that $K_s\cap Z_{2i-s+1}(G)$ has finite index in $K_s$. For $s=i+1$ the
statement is true, as $|K_{i+1}:K_{i+1}\cap Z_i(G)|\le
|\gamma_{i+1}(G):\gamma_{i+1}(G)\cap Z_i(G)|$ is finite by hypothesis. Now assume
that $Z=K_{s+1}\cap Z_{2i-s}(G)$ has finite index in $K_{s+1}$.
As $C$ has finite index in $G$, we have $G=\langle g_1,\ldots,g_n,C\rangle$
for some $g_1,\ldots,g_n\in G$. Let $U=\langle
g_1,\ldots,g_n\rangle$ and let
$H_s$ be the centralizer of $UZ/Z$ in $K_s$. Since
$[K_s,U]Z/Z\le K_{s+1}/Z$ is finite and $U$ is finitely generated,
the subgroup $H_s$ has finite index in $K_s$. Moreover, $[H_s,U]\le Z\le Z_{2i-s}(G)$ and
$[H_s,C]\le [K_s,C]\le [C,\,_{s-1}G,C]\le Z_{2i-s}(G)$, where the last inclusion follows
from Lemma \ref{comm}.
Hence $[H_s,G]=[H_s,UC]\le Z_{2i-s}(G)$ and $H_s\le Z_{2i-s+1}(G)$, which completes the induction.

In particular, for $s=1$ we obtain that $|C:C\cap Z_{2i}(G)|$ is finite.
Consequently $|G:Z_{2i}(G)|\le |G:C\cap Z_{2i}(G)|=|G:C|\,|C:C\cap Z_{2i}(G)|$
is finite, and we are done.
\end{proof}
\bibliographystyle{amsplain}

\end{document}